\def\mytitle{A Survey on Dual-Quaternions} 

\def\mysubtitle{}

\def\mykeywords{survey, sound, research, kinematics, computer graphics, dual quaternions, complex numbers, trends, publications, impact}

\def\myauthor{Benjamin Kenwright}
\def\myemail{ }

\documentclass[9pt,journal,cspaper,compsoc]{IEEEtran}

\usepackage{times}

\usepackage{graphicx}
\usepackage{watermark}

\usepackage[backref=page]{hyperref}

\hypersetup{
    colorlinks = true, 
    linkcolor = blue,
    anchorcolor = red,
    citecolor = blue, 
    filecolor = red, 
}

\usepackage{tikz}
\usepackage{calc}
\usepackage{setspace}

\usetikzlibrary{math}

\usepackage{times}

\newcommand{\drawledger}[1]
{
	\draw (0+#1,0) -- (0+#1,20.0cm);
	
	\foreach \x in {0,1,...,20} 
	{
		\draw (3pt+#1, \x cm) -- (-3pt+#1, \x cm);
		
	}
	
	\draw (#1,0cm) node[left=3pt,font=\fontsize{7pt}{7pt}\selectfont] {$ 2010 $} ;
	\draw (#1,20cm) node[left=3pt,font=\fontsize{7pt}{7pt}\selectfont] {$ 1950 $} ;

}

\usepackage{tikz}
\usetikzlibrary{shapes,arrows}

\usepackage{supertabular}
\usepackage{multicol}

\usepackage{tikz}
\usepackage{calc}
\usepackage{setspace}

\usetikzlibrary{mindmap}

\usepackage{pgfplots}

\usepackage{pgf-pie}

\usepackage{smartdiagram}

\usepackage{bchart}

\usepackage{subcaption}

\usepackage{graphicx}
\graphicspath{{./images/}}

\newcommand{\figuremacroW}[4]{
	\begin{figure}[htbp]
		\centering
		\includegraphics[width=#4\columnwidth]{#1}
		\caption[#2]{\textbf{#2} - #3}
		\label{fig:#1}
	\end{figure}
}

\newcommand{\figuremacroF}[4]{
	\begin{figure*}[htbp]
		\centering
		\includegraphics[width=#4\textwidth]{#1}
		\caption[#2]{\textbf{#2} - #3}
		\label{fig:#1}
	\end{figure*}
}

\usepackage{color}  
\usepackage{xcolor}
\definecolor{lbcolor}{rgb}{0.98,0.98,0.98}

\usepackage{listings}

\newcommand{\subparagraph}{}

\usepackage{titlesec}

\renewcommand\theparagraph{}

\titleformat*{\paragraph}{\bfseries}
\titleformat{\paragraph}[runin]
	{\normalfont\normalsize\bfseries}
	{\theparagraph}
	{0em}
	{}

\titlespacing*{\paragraph}
{0pt}
{0.25ex plus 1ex minus .2ex}
{1em}

\usepackage{tocloft}

\usepackage{amsmath} 

\usepackage{amssymb}
\usepackage{amsmath}

\begin{document}

\raggedbottom
\sloppy

\title{\fontsize{16}{22}\selectfont \mytitle \\ \fontsize{13}{26}\selectfont \mysubtitle}



\hypersetup{pdfinfo={
   Author		= {\myauthor},
   Title		= {\mytitle  \mysubtitle},
   Subject 		= {\mytitle  \mysubtitle},
   CreationDate = {D:20230220195600},
   Keywords 	= {\mykeywords},
}}

\author{\myauthor
\IEEEcompsocitemizethanks{\IEEEcompsocthanksitem \myauthor \; (2023) \protect\\
\myemail 
}
\thanks{}}

\raggedbottom

\markboth{Communication Article (\myauthor) - \mytitle \; (2023)}%
{Technical Paper}

\IEEEcompsoctitleabstractindextext{%
\begin{abstract}
\boldmath
%
Over the past few years, the applications of dual-quaternions have not only developed in many different directions but has also evolved in exciting ways in several areas.
%
%
As dual-quaternions offer an efficient and compact symbolic form with unique mathematical properties. 
%
While dual-quaternions are now common place in many aspects of research and implementation, such as, robotics and engineering through to computer graphics and animation, there are still a large number of avenues for exploration with huge potential benefits.
This article is the first to provide a comprehensive review of the dual-quaternion landscape. 
In this survey, we present a review of dual-quaternion techniques and applications developed over the years while providing insights into current and future directions. 
The article starts with the definition of dual-quaternions, their mathematical formulation, while explaining key aspects of importance (e.g., compression and ambiguities). 
The literature review in this article is divided into categories to help manage and visualize the application of dual-quaternions for solving specific problems. 
A timeline illustrating key methods is presented, explaining how dual-quaternion approaches have progressed over the years. 
The most popular dual-quaternion methods are discussed with regard to their impact in the literature, performance, computational cost and their real-world results (compared to associated models).
Finally, we indicate the limitations of dual-quaternion methodologies and propose future research directions.

\end{abstract}
\begin{keywords}
\mykeywords
\end{keywords}}

\maketitle

\IEEEdisplaynotcompsoctitleabstractindextext
\IEEEpeerreviewmaketitle

\section{Introduction}

\paragraph{Evolving Technology}

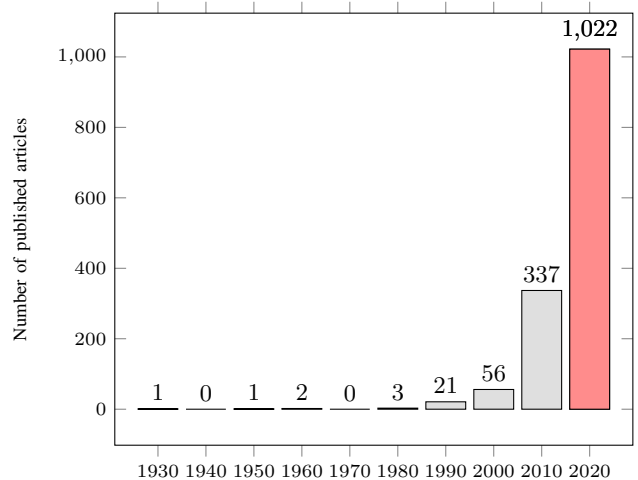
\begin{figure}
\begin{tikzpicture}
\begin{axis}[
    ybar,
    ylabel={Number of published articles},
    ylabel style={font={\scriptsize}},
    xtick=data,
    xticklabel style={/pgf/number format/1000 sep=,{font=\scriptsize}},
    yticklabel style={font=\scriptsize},
    nodes near coords,
    nodes near coords align={vertical},
    ybar=0pt, bar width=15, bar shift=0pt
    ]
\addplot [fill=gray!25] coordinates 
{   
    (1930 , 1)
    (1940 , 0)
    (1950 , 1)
    (1960 , 2)
    (1970 , 0)
    (1980 , 3)
    (1990 , 21)
    (2000 , 56)
    (2010 , 337)
    (2020 , 1022)
};
\addplot [fill=red!45] coordinates 
{   
    (2020 , 1022)
};
\end{axis}
\end{tikzpicture}
\caption{Visually illustrating the rising popularity and impact of dual-quaternions by showing the number of articles with the word 'dual-quaternion' in the title published between 1930 and 2020 (predictions upto 2030). (Google scholar data accessed 23/03/2023 keywords 'dual-quaternion' in title). }
\label{fig:numpapers}
\end{figure}

The technological advancements in quantum computing \cite{bravyi2022future}, 3D-printing \cite{zhu2022recent}, flexible transparent screens \cite{wang2021recent}, and breakthroughs in machine learning and artificial intelligence \cite{ray2019quick} offer new possibilities and problems while requiring new approaches to be developed and investigated.
One such approach is dual-quaternions.
Dual-quaternions are a combination of dual-number theory and complex-numbers which offer a novel mathematical form with a unique set of properties.
While dual-quaternion's have made a substantial impact in robotics over the years, they have recently found applications in other fields, such as in computer graphics with particular interest in the area of animation. 
This survey explores the dual-quaternion landscape, surveying some of the most popular applications, problems and trends.
The purpose is to explain and demonstrate the diverse abilities of dual-quaternions to solve a range of problems, in some case having a substantial impact. Over the past decade, there has been a substantial increase in papers with dual-quaternions in the title, yet there has been to date no comprehensive surveys on dual-quaternions (i.e., review of the topic, properties, limitations, applications and comparisons of dual-quaternion techniques holistically).
Of the 700+ publications with dual-quaternions in the title since 1930, the majority of them are from the last few years (Figure \ref{fig:numpapers}); showing the sudden explosion of novel research around dual-quaternions in multiple areas like sound and graphics. 
As we show in this article, the multidisciplinary value of dual-quaternions and its initial benefits in kinematics and robotics have been a catalyse for its continued growth and success in many domains.

\paragraph{Contribution}
The contribution of this article is a comprehensive review of dual-quaternion techniques, which is the first survey on dual-quaternions, including recent trends and new future directions.
Dual-quaternion techniques are presented including the categorization and application of them.
We highlights the advantages and disadvantages of dual-quaternion approaches with regard to comparable methods in context (e.g., robustness, complexity and computational costs).
We also provide indications on which dual-quaternion techniques are the best fit and why for solving specific problems, while discussing opportunities for future work and applications (based on similar research and/or models). 

\section{Background Mathematics}
To ensure this article provides a self-contained review of the topic, we briefly describe they key mathematical concepts behind dual-quaternions.
Dual-quaternions, although not as well known as Quaternions, provide a fundamental and solid base to describing three-dimensional transform (orientation and translation) of an object or a vector. 
They are efficient and well suited to solve a variety of problems in computer graphics and animation.

\begin{figure*}
    \centering
    \tikzset{every picture/.style={line width=0.75pt}} 

\begin{tikzpicture}[x=0.75pt,y=0.75pt,yscale=-1,xscale=1]

\draw  [color={rgb, 255:red, 74; green, 144; blue, 226 }  ,draw opacity=1 ][dash pattern={on 4.5pt off 4.5pt}] (100,196.2) .. controls (100,186.7) and (107.7,179) .. (117.2,179) -- (523.59,179) .. controls (533.09,179) and (540.79,186.7) .. (540.79,196.2) -- (540.79,247.8) .. controls (540.79,257.3) and (533.09,265) .. (523.59,265) -- (117.2,265) .. controls (107.7,265) and (100,257.3) .. (100,247.8) -- cycle ;
\draw  [color={rgb, 255:red, 74; green, 144; blue, 226 }  ,draw opacity=1 ][dash pattern={on 4.5pt off 4.5pt}] (100,95) .. controls (100,85.06) and (108.06,77) .. (118,77) -- (522.79,77) .. controls (532.73,77) and (540.79,85.06) .. (540.79,95) -- (540.79,149) .. controls (540.79,158.94) and (532.73,167) .. (522.79,167) -- (118,167) .. controls (108.06,167) and (100,158.94) .. (100,149) -- cycle ;
\draw  [color={rgb, 255:red, 128; green, 128; blue, 128 }  ,draw opacity=1 ] (260.19,118.2) .. controls (260.19,116.29) and (261.75,114.73) .. (263.66,114.73) -- (274.06,114.73) .. controls (275.97,114.73) and (277.53,116.29) .. (277.53,118.2) -- (277.53,129.33) .. controls (277.53,131.25) and (275.97,132.8) .. (274.06,132.8) -- (263.66,132.8) .. controls (261.75,132.8) and (260.19,131.25) .. (260.19,129.33) -- cycle ;
\draw  [color={rgb, 255:red, 128; green, 128; blue, 128 }  ,draw opacity=1 ] (288.59,118.35) .. controls (288.59,116.35) and (290.21,114.73) .. (292.21,114.73) -- (362.84,114.73) .. controls (364.83,114.73) and (366.45,116.35) .. (366.45,118.35) -- (366.45,129.19) .. controls (366.45,131.18) and (364.83,132.8) .. (362.84,132.8) -- (292.21,132.8) .. controls (290.21,132.8) and (288.59,131.18) .. (288.59,129.19) -- cycle ;
\draw  [color={rgb, 255:red, 128; green, 128; blue, 128 }  ,draw opacity=1 ] (312.59,219.85) .. controls (312.59,217.85) and (314.21,216.23) .. (316.21,216.23) -- (336.7,216.23) .. controls (338.7,216.23) and (340.32,217.85) .. (340.32,219.85) -- (340.32,230.69) .. controls (340.32,232.68) and (338.7,234.3) .. (336.7,234.3) -- (316.21,234.3) .. controls (314.21,234.3) and (312.59,232.68) .. (312.59,230.69) -- cycle ;
\draw  [color={rgb, 255:red, 128; green, 128; blue, 128 }  ,draw opacity=1 ] (277.82,219.85) .. controls (277.82,217.85) and (279.44,216.23) .. (281.43,216.23) -- (300.2,216.23) .. controls (302.2,216.23) and (303.82,217.85) .. (303.82,219.85) -- (303.82,230.69) .. controls (303.82,232.68) and (302.2,234.3) .. (300.2,234.3) -- (281.43,234.3) .. controls (279.44,234.3) and (277.82,232.68) .. (277.82,230.69) -- cycle ;
\draw    (320.4,225.05) -- (324.28,245.5) ;

\draw    (337.19,200.55) -- (327.69,215.59) ;

\draw    (290.19,203.55) -- (292.69,216.09) ;

\draw    (349.9,129.05) -- (356.79,148.11) ;

\draw    (342.19,99.55) -- (332.69,114.59) ;

\draw    (265.56,102.19) -- (268.06,114.73) ;

\draw (299.33,123) node  [align=left] {q = w + \textit{\textcolor[rgb]{0.29,0.56,0.89}{i}}x + \textit{\textcolor[rgb]{0.29,0.56,0.89}{j}}y + \textit{\textcolor[rgb]{0.29,0.56,0.89}{k}}z};
\draw (296.82,223.27) node  [align=left] {$\displaystyle \zeta =q_{r} \ +\epsilon q_{d}$};
\draw (147,91) node  [align=left] {Quaternion};
\draw (168,197) node  [align=left] {Dual-Quaternion};
\draw (479,251) node [scale=0.9] [align=left] {8 scalar variables};
\draw (477,155) node [scale=0.9] [align=left] {4 scalar variables};
\draw (262,95) node [scale=0.9] [align=left] {Real};
\draw (353.19,193.55) node [scale=0.9] [align=left] {Dual-Part};
\draw (287,194) node [scale=0.9] [align=left] {Real};
\draw (344.19,94.55) node [scale=0.9] [align=left] {Complex};
\draw (360.19,152.55) node [scale=0.8] [align=left] {Imaginary};
\draw (336.19,253.05) node [scale=0.8] [align=left] {Dual Operator};

\end{tikzpicture}
    \caption{Visual Overview of Quaternion and Dual-Quaternion Components.}
    \label{fig:dualquatlayout}
\end{figure*}
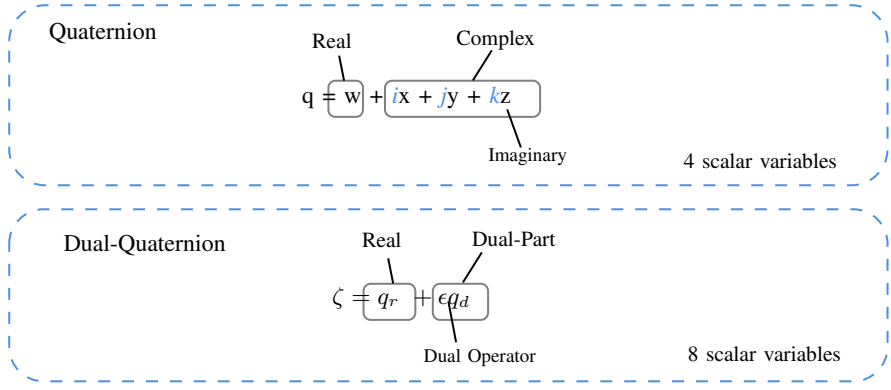

\subsection{Algebraic Definitions}

\paragraph{Quaternion Operations}

The quaternions were discovered by Hamilton in 1843 as a method of performing 3-D multiplication \cite{kenwright2012beginners}. 
A quaternion $q$ is given by Equation \ref{eq:quatdef}.
Since we are combining quaternions with dual number theory, we give the elementary quaternion arithmetic operations.

\begin{equation}
q = [s,\vec{v}], (q_w,q_x,q_y,q_z)
\label{eq:quatdef}
\end{equation}

where $s$ scalar part is $s=q_w$ and vector part is $\vec{v}=(q_x,q_y,q_z)$. The four-tuple of independent real values assigned to one real axis and three orthonormal imaginary axes: $i,j,k$.

\begin{itemize}

\item \textbf{addition}: $q_1 + q_2 = [s_1,\vec{v}_1]+[s_2,\vec{v}_2] = [s_1+s_2,\vec{v}_1+\vec{v}_2]$

\item \textbf{additive identity}: $0 = [0,0]$

\item \textbf{scalar multiplication}: $kq=[ks,k\vec{v}]$

\item \textbf{multiplication}: $q_1 q_2 = [s_1,\vec{v}_1][s_2,\vec{v}_2] = [s_1s_2 - \vec{v}_1 \cdot \vec{v}_2, s_1\vec{v}_2 + s_2\vec{v}_1 + \vec{v}_1 \times \vec{v}_2]$

\item \textbf{multiplication identity}: $1 = [1,0]$

\item \textbf{dot product}: $q_1 \cdot q_2 = ( q_{1x}q_{2x} + q_{1y}q_{2y} + q_{1z}q_{2z} + q_{1w}q_{2w} )$

\item \textbf{magnitude}: $||q|| = \sqrt( s^2 + ||\vec{v}||^2 )$

\item \textbf{conjugate} $q^* = [s,-\vec{v}]$

\end{itemize}

\paragraph{Dual-Quaternion Operations}

The elementary arithmetic operations necessary for us to use dual-quaternions.

\begin{itemize}

\item \textbf{dual-quaternion}: $\zeta = q_r + q_d \varepsilon$

\item \textbf{scalar multiplication}:  $s \zeta = s q_r + s q_d \varepsilon$

\item \textbf{addition}: $\zeta_1 + \zeta_2 = q_{r1} + q_{r1} + (q_{d1}+q_{d2}) \varepsilon$

\item \textbf{multiplication}: $\zeta_1 \zeta_2 = q_{r1} + q_{r2} + (q_{r1}q_{d2} + q_{d1}q_{r2}) \varepsilon$

\item \textbf{conjugate}:  $\zeta^* = q_r^* + q_d^* \varepsilon$

\item \textbf{magnitude}: $||\zeta|| = \zeta \zeta*$

\end{itemize}
\noindent where $q_r$ and $q_d$ indicate the real and dual part of a dual-quaternion.

For a beginners introduction to dual-quaternions and a comparison of alternative methods (e.g., matrices and Euler angles) and how to go about implementing a straightforward library we refer the reader to the paper by Kenwright \cite{kenwright2012beginners}.

\subsection{Algebra in Context}

\paragraph{Dual-Quaternion Vector Transformation}
A dual-quaternion is able to transform a 3D vector coordinate as shown in Equation 34. Note that for a unit-quaternion the inverse is the same as the
conjugate.

\begin{equation}
    p' = \hat{\zeta} p \hat{\zeta}^{-1}
\end{equation}

\noindent where $\hat{\zeta}$ is a unit dual-quaternion representing the transform, $\hat{\zeta}^{-1}$ is the inverse of the unit dual-quaternion transform.
$p$ and $p'$ are the dual-quaternions holding 3D vector coordinate to before and after the transformation (i.e., $p=(1,0,0,0) + \epsilon (0,v_x,v_y,v_z)$ )).

\paragraph{Pl{\"u}cker Coordinates}
Pl{\"u}cker coordinates are used to create Screw coordinates which are an essential technique of representing lines. We need the Screw coordinates so that we can re-write dual-quaternions in a more elegant form to aid us in formulating a neater and less complex interpolation method that is comparable with spherical linear interpolation for classical quaternions.

The Definition of Plücker Coordinates:
\begin{itemize}
\item $\vec{p}$ is a point anywhere on a given line
\item $\vec{l}$ is the direction vector
\item $\vec{m} = \vec{p} \times \vec{l}$ is the moment vector
\item $(\vec{l},\vec{m})$ are the six Plücker coordinate
\end{itemize}

We can convert the eight dual-quaternions parameters to
an equivalent set of eight screw coordinates and vice-versa.
The definition of the parameters are given in Equation

\begin{equation}
\begin{alignedat}{4}
\text{screw parameters} &= (\theta,d,\vec{l},vec{m}) \\
\text{dual-quaternion}  &= q_r + \epsilon q_d \\
                        &= (w_r+\vec{v}_r) + \epsilon(w_d + \vec{v}_d) 
\end{alignedat}
\end{equation}

\noindent where in addition to $\vec{l}$ representing the vector line direction and $\vec{m}$ the line moment, we also have $d$
representing the translation along the axis (i.e., pitch) and the angle of rotation $\theta$.

\textbf{Convert dual-quaternion to screw-parameters}

\begin{equation}
\begin{alignedat}{3}
\theta &= 2 cos^{-1}(w_r) \\
d      &= -2w_d \frac{1}{\sqrt(\vec{v}_r \dot \vec{v}_r)} \\
\vec{l} &= \vec{v}_r \left( \frac{1}{\sqrt(\vec{v}_r \dot \vec{v}_r)} \right) \\
\vec{m} &= \left( \vec{v}_d - \vec{l} \frac{d w_r}{2} \right) \frac{1}{\sqrt(\vec{v}_r \dot \vec{v}_r)}
\end{alignedat}
\end{equation}

\textbf{Convert screw-parameters to dual-quaternion}

\begin{equation}
\begin{alignedat}{3}
w_r &= cos \left( \frac{\theta}{2} \right) \\
\vec{v}_r &= \vec{l} sin \left( \frac{\theta}{2} \right) \\
w_d &= -\frac{d}{2} sin \left( \frac{\theta}{2} \right) \\
\vec{v}_d &= sin \left( \frac{\theta}{2} \right) \vec{m} + \frac{d}{2} cos \left( \frac{\theta}{2} \right) \vec{l}
\end{alignedat}
\end{equation}

\paragraph{Dual-Quaternion Power}

We can write the dual-quaternion representation in the
form given in Equation \ref{eq:dqpow}. 

\begin{equation}
\begin{alignedat}{3}
\hat{\zeta} &= cos \left( \frac{\theta + \epsilon d}{2} \right)
            + (\vec{l} + \epsilon \vec{m}) 
               sin \left( \frac{\theta + \epsilon d}{2} \right) \\
            &= cos \left( \frac{\hat{\theta}}{2} \right) +
                \hat{v} sin \left( \frac{\hat{\theta}}{2} \right)
\end{alignedat}
\label{eq:dqpow}
\end{equation}

\noindent where $\hat{\zeta}$ is a unit dual-quaternion, $\hat{v}$ is a unit dual-vector ($\hat{v}=\vec{l}+\epsilon\vec{m}$), and $\hat{\theta}$ is a dual-angle ($\hat{\theta}=\theta+\epsilon d$).

The dual-quaternion in this form is exceptionally
interesting and valuable as it allows us to calculate a dual-quaternion to a power. Calculating a dual-quaternion to a
power is essential for us to be able to easily calculate
spherical linear interpolation. However, instead of purely
rotation as with classical quaternions, we are instead now
able to interpolate full rigid transformations (i.e., rotation
and translation) by using dual-quaternions. 

\begin{equation}
    \hat{\zeta}^{t} = cos \left( t \frac{\hat{\theta}}{2} \right) + 
    \hat{v} sin \left( t \frac{\hat{\theta}}{2} \right)
\end{equation}

\paragraph{Dual-Quaternion Screw Linear Interpolation (ScLERP)}
ScLERP is an extension of the quaternion SLERP
technique, and allows us to create constant smooth
interpolation between dual-quaternions. 
Similar to quaternion SLERP we use the power function to calculate
the interpolation values for ScLERP shown in Equation

\begin{equation}
\begin{alignedat}{3}
\text{ScLERP}(\hat{\zeta}_A, \hat{\zeta}_B : t ) = \hat{\zeta}_A (\hat{\zeta}^{-1}_A \hat{\zeta}_B )^t
\end{alignedat}
\end{equation}

\noindent where
$\hat{\zeta}_A$
and
$\hat{\zeta}_B$
are the start and end unit dual-quaternion and $t$ is the interpolation amount from $0.0$ to $1.0$.

Alternatively, a fast approximate alternative to ScLERP was presented by Kavan et al. \cite{pop00001} called Dual-Quaternion Linear Blending (DLB). Furthermore, dual-quaternions
have gained a great deal of attention in the area of character-based skinning. Since, a skinned surface
approximation using a weighted dual-quaternion approach produces less kinking and reduced visual anomalies
compared to linear methods by ensuring the surface keeps its volume.

There have been elaborations of the concept for artistic control to allow blending, amplifying or mixing aspects of the interpolation as presented by Kenwright \cite{kenwright2023dualquaternion} (e.g., KenLERP).


\figuremacroW
{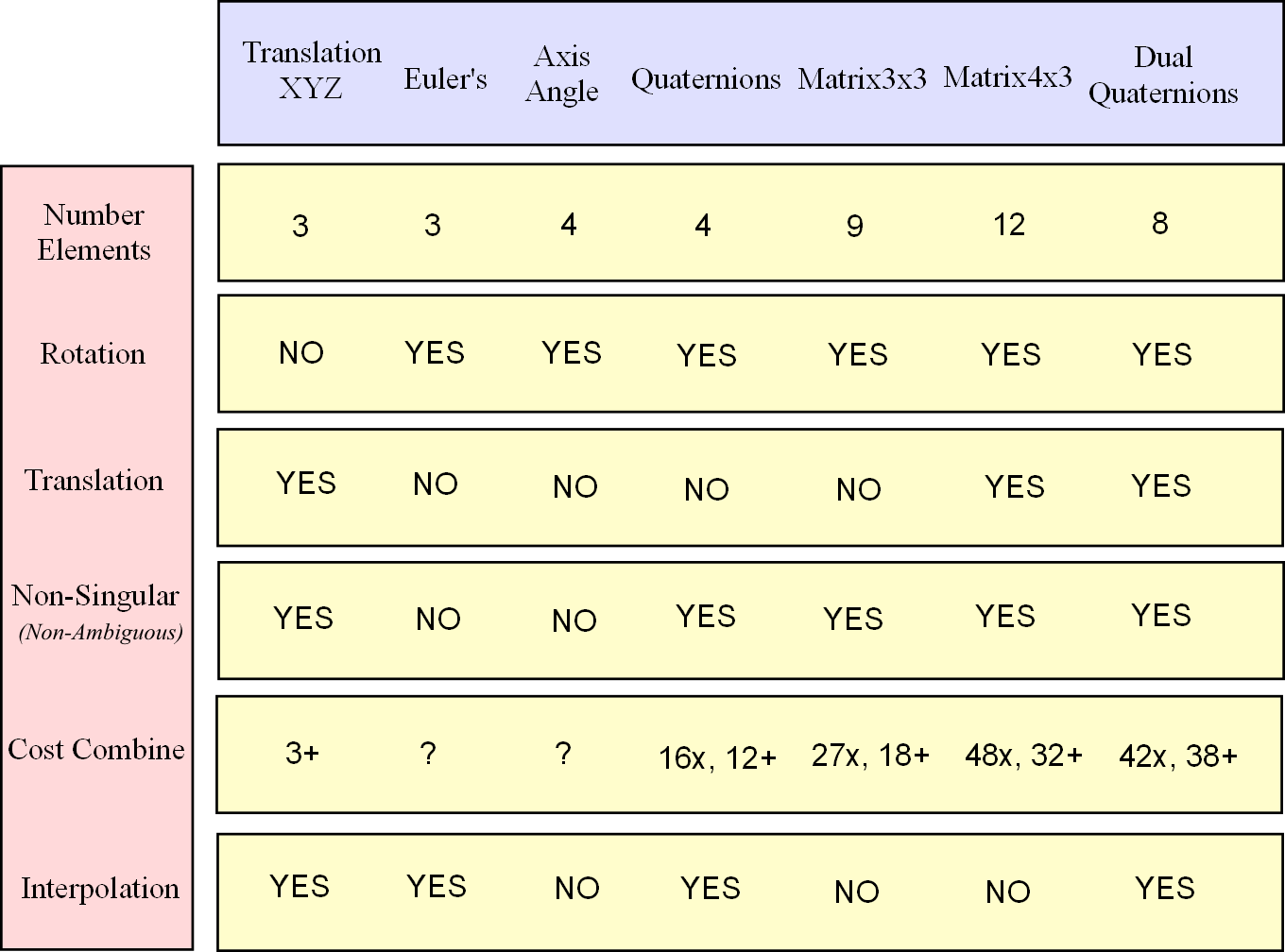}
{Geometric Transformation}
{Comparison of Similar Methods Transformation \cite{kenwright2012beginners}. Note the $+/\times$ symbols in the table for the computational cost refer to the number of multiplications and additions.}
{1.0}

\subsection{Computational Factors}
The computational cost of dual-quaternions is related to a number of factors, such as, the hardware acceleration model. While Table \ref{fig:dqtable} provides a `arithmetic' cost for combining objects, it does not take into account the other operations, such as, interpolation or manipulations that influence the overall 'real-world' computational cost. The implementation and data may have influence run-time testing (caching, loops and how the operations are performed).



\section{Applications}

\subsection{Dual-Quaternion Applications in Computer Vision}
In computer vision, dual-quaternions have been used to solve problems such as registration of multiple views \cite{dod00011} and hand-eye calibration \cite{dod00001}.

\paragraph{Hand-Eye Calibration}
Daniilidis \cite{dod00001} related measurements from sensors on physical robots to different coordinate frames (e.g., point on a 2d image to a 3d location). 
Dual-quaternions were able to represent this hand-eye calibration transform in a unified manner.

\paragraph{Mult-View Reconstruction}
Torsello et al. \cite{dod00011} proposed a novel multi-view reconstruction algorithm using dual-quaternions for reconstructing three-dimensional surfaces from 2-dimensional image views. 
The approach allowed for a completely generic topology and through extensive experimentation was shown to be both orders of magnitude faster and more robust to extreme positional noise than other approaches. 

\figuremacroF
{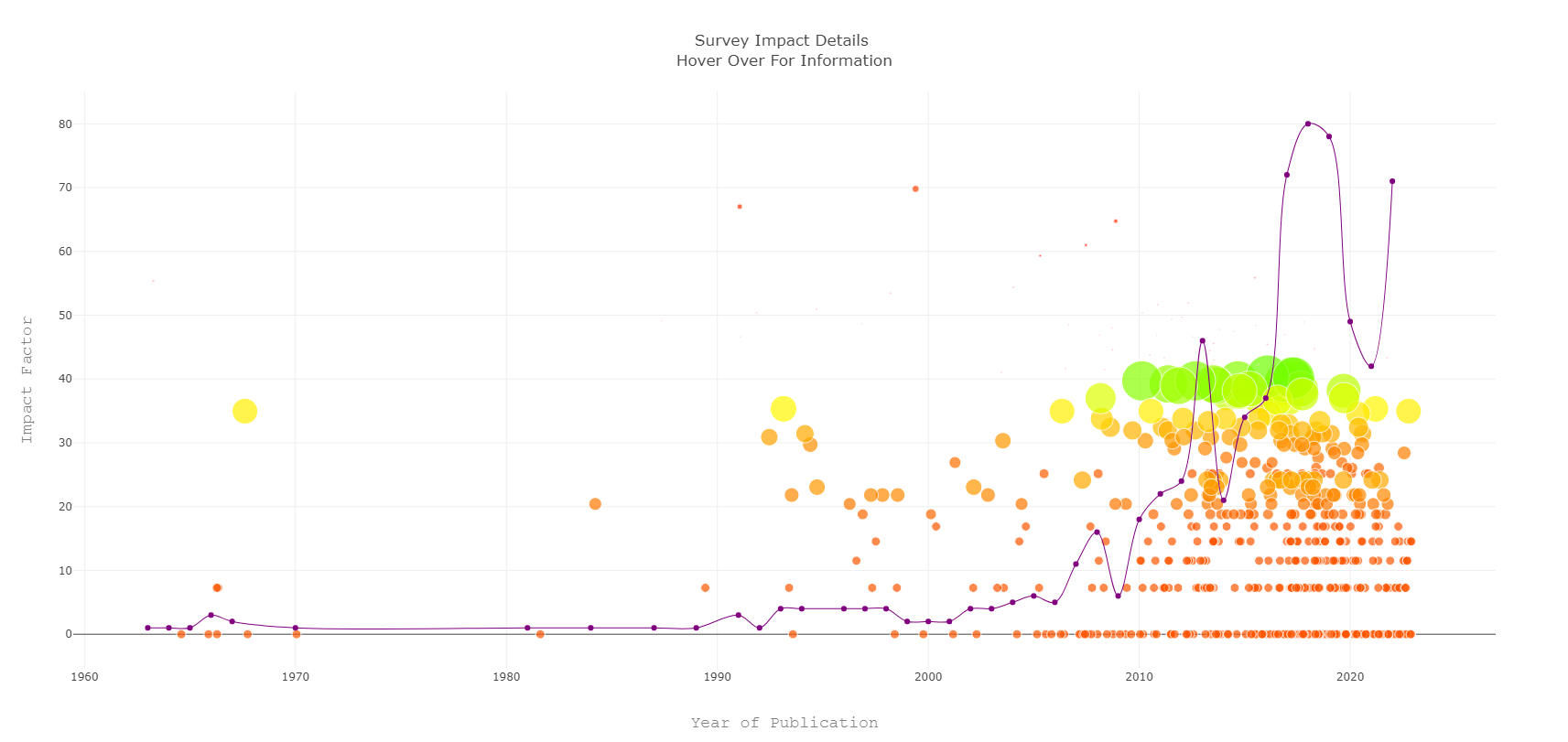}
{Timeline}
{Search terms: dual-quaternion or dual-quaternions or 'dual quaternion' or 'dual quaternions'.
	No Articles: 709.
	Retrieved Date: 23/03/2023 \cite{dqdata2023}. }
{1.0}

\subsection{Dual-Quaternion Applications in Kinematics}

The conventional method for representing and concatenating links together in hierarchical systems is the Denavit-Hartenberg matrix convention, and while Wang and Ravani proposed an alternative more efficient forward recursion method for kinematic equations, we propose using dual-quaternions, since they offer an analogous alternative that is numerically stable and computationally efficient. Dual-quaternions have shown promising results for providing singularity-free solutions for inverse kinematic (IK) problems with non-linearities. It is clearly an advantage to use dual-quaternions for rigid hierarchies since each dual-quaternion can be concatenated easily, interpolated smoothly and provide rigid transform comparisons effortlessly. 
Generating fast reliable Inverse Kinematic (IK) solutions in real-time with angular limits for highly articulated figures (e.g., human bipeds including hands and feet) is challenging and important. The subject is studied across numerous disciplines, such as graphics, robotics, and biomechanics, and is employed by numerous applications in entertainment, training and medical fields \cite{Kenwright2013InverseKW}.

\subsubsection{Dual-Quaternions Achilles' Heel (Limitations)}
Dual-quaternions have limitations, when used as geometric `transformations', they are unable to represent non-rigid concepts (the `rotational' part is assumed to be a unit quaternion, the translational part is essentially a translation vector - so there is no room for skew and scale).
This means that for extrinsic parameters, dual quaternions are a good representation, while intrinsics are problematic.

\subsubsection{Data Binding and Representation}
The majority of research and applications for dual-quaternions have been credited to geometric transformations (compact method for coupling rotational and translational data into a single form that can easily be combined and manipulated/interpolated). However, crucially, the value and application of dual-quaternions is not limited by these data formats. 
For instance, dual-quaternions can be used to represent `color' information (red, green, blue) (combined with displacement and/or pattern characteristics). Alternatively it can be used for describing `signals' such as those used in sound and audio.
The underlying mathematics and theorems can also be blended and combined with analogues processes to reap additional benefits (Fourier and Neural Networks).

\subsection{Dual-Quaternion Applications in Graphics and Other Areas}

Examples of applications of dual-quaternions both directly and indirectly based on studies, research or similarity with models and/or mathematical principles:

\begin{itemize}
	\item Sound and audio related processes \cite{kenwright2023dualquaternions,grassucci2023dual}
    \item Image processing (Dual-quaternion Fourier \cite{ding2012approach}) 
    \item Character Skinning \cite{pop00001}
    \item Articulated Skeletons \cite{kenwright2012beginners}
    \item Kinematics and Inverse-Kinematics \cite{Kenwright2013InverseKW}
    \item Physics (Classical Mechanics) \cite{goddard1998pose}
    \item Filtering/Compression \cite{sveier2020dual}
    \item Interpolation \cite{kenwright2023dualquaternion}
    \item Geometric (Multiview Reconstruction \cite{dod00011})
    \item Fractals (Visualization Julia \cite{kenwrightdualfractals})
    \item Optimization, transmission/storage, noise 
    \item Spline Curves and Surfaces (NURBS) \cite{kenwright2018dual}
    \item Machine learning (neural networks) \cite{kenwright2023dualquaternions}
\end{itemize}

\section{Taking Dual-Quaternions to the Edge (Future Directions)}
While dual-quaternions have made an impact in multiple areas, like kinematics, animation and robotics, there are still unexplored areas for potential opportunities, such as, security, image processing, compression and computational-models.

While this article has highlighted the incredible innovations and advancements in the literature using dual-quaternions, it only scratches the surface. Dual-quaternions are a flexible tool that can be used in a vast assortment of ways.  For example, dual-quaternions offer an alternative solution for representing, manipulating or combining data in a unique way to solve problems or provide enhancements (these benefits are far and flexible - from computational to analytical complexity factors).


\begin{figure*}
	\centering
	\input{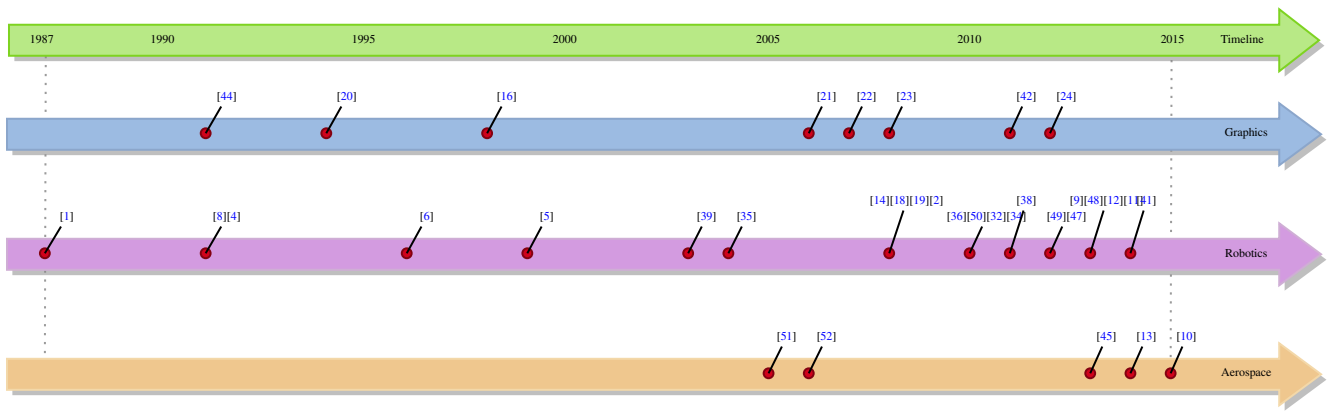}
	\caption{Timeline of key articles published over the past few decades that have made a significant impact on the topic of dual-quaternions (in graphics, robotics and aerospace). }
	\label{fig:mytimeline}
\end{figure*}


\subsection{Security}
While security has not been touched on much with dual-quaternions - there is potential opportunities, especially when looking at related research that could be extended from quaternion-based models.

For example, signals such as audio and animation, transformed into dual-quaternion space may offer a novel Steganography solution (tagging or tracking data). One possible way of accomplishing this would be to extend the quaternion-based research which used geometric invariance of the quaternion exponent to embedding secret information \cite{li2021image}.

Another avenue of exploration maybe, dual-quaternion public key encryption schemes, using similar work that follows on from quaternion key models \cite{valluri2016quaternion}, could be extended to dual-quaternions; allowing the based idea to be elaborated and enhanced. 

\figuremacroW
{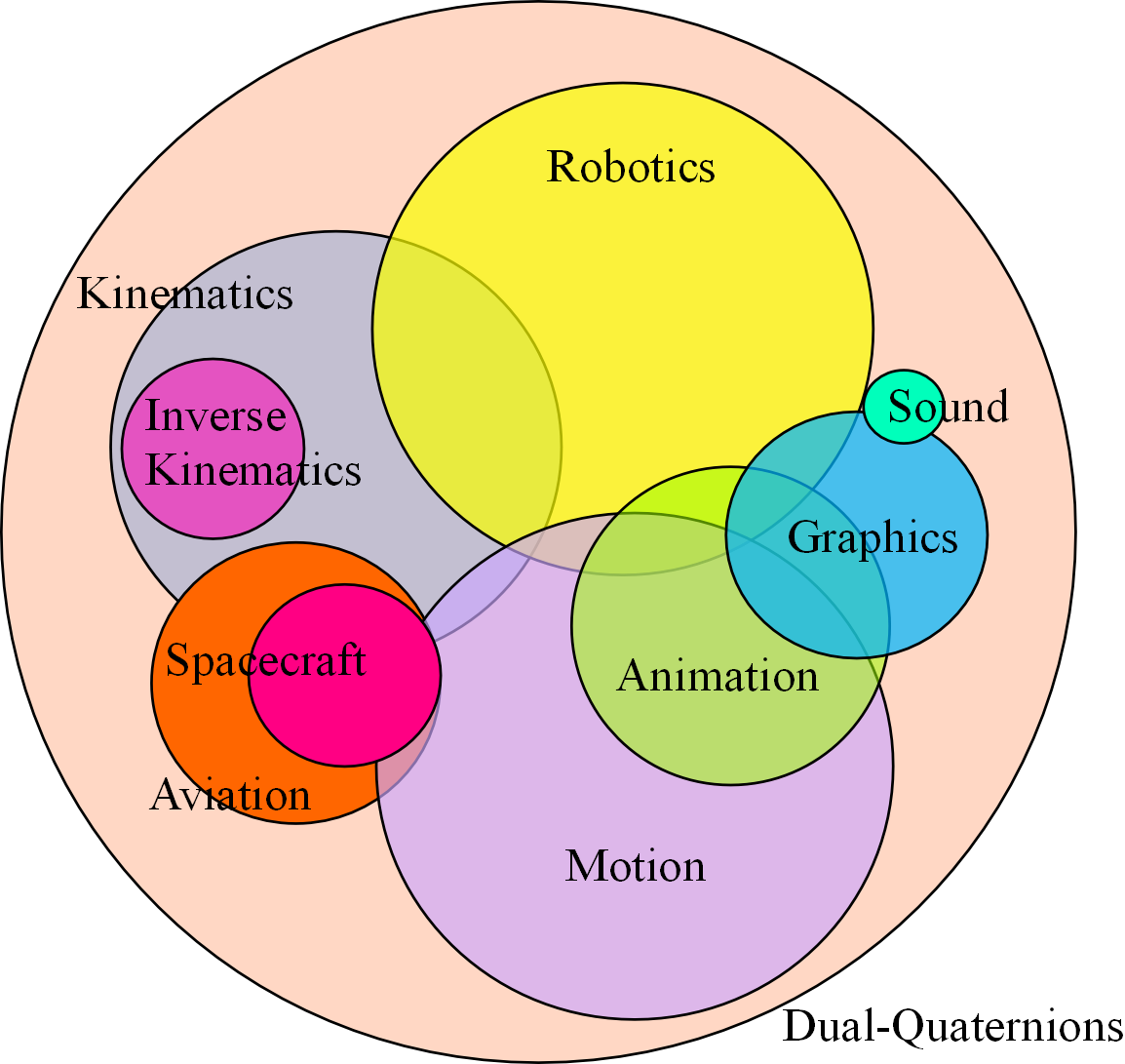}
{Topic Analysis}
{Topic model analysis using the title, abstract text and keywords, we group and visualize key contributions to date that have used dual-quaternions.  }
{0.9}

\section{Conclusion and Discussion}

In this paper, we have surveyed some of the most popular dual-quaternion techniques, emphasizing, but not limiting, the discussion to approaches from a computer graphics point of view. 
We describe where the research in dual-quaternions has been focused in the past, how it has progressed over the years and indicate reasons for such progression. 
The main scope of this survey is to offer a guide that highlights the advantages and disadvantages of important dual-quaternion techniques, giving indications about which methods are best suited to solve specific problems. 
It aims to introduce dual-quaternions to new researchers that aim to optimize their projects. 
Finally, it provides future directions to extend current limitations, and research challenges that need further investigation.


\bibliographystyle{plain}





%
%

\bibliography{paper} 


\end{document}